\documentclass{amsproc}

\newtheorem{theorem}{Theorem}[section]

\theoremstyle{definition}

\newtheorem{example}[theorem]{Example}
\theoremstyle{remark}
\newtheorem{remark}[theorem]{Remark}
\usepackage{rotating}

\numberwithin{equation}{section}

\usepackage{mathrsfs} 
\usepackage{graphicx} 
\usepackage{graphicx,graphics,epsf,psfrag,eepic,epic} 

\def\N{\mathbb{Z}_+}
\def\Z{\mathbb{Z}}
\def\R{\mathbb{R}}
\def\C{\mathbb{C}}
\def\g{\mathfrak{g}}
\def\t{\mathfrak{t}}
\def\la{\lambda}

\def\eps{\varepsilon}
\def\res{\operatorname{Res}}
\def\ires{\operatorname{IRes}}

\def\spa{\operatorname{Sp}({\bf a})}
\def\bnum{\begin{enumerate}}
\def\enum{\end{enumerate}}
\def\clm{c_\la^\mu}
\def\ctlm{c_{t\la}^{t\mu}}
\def\clmn{c_{\la\,\,\mu}^{\ \, \nu}}
\def\ctlmn{c_{t\la\ \,t\mu}^{\ \ t\nu}}
 
\def\klr{weight multiplicities and Littlewood-Richardson coefficients} 
\def\ch{\operatorname{ch}}
\def\slr#1{\mathfrak{sl}_{#1}(\C)}
\def\glr#1{\mathfrak{gl}_{#1}(\C)}

\def\vlm{\operatorname{Val}(\la,\mu)}
\def\vlmn{\operatorname{Val}(\la,\mu,\nu)}
\def\ppa{P(\Phi,a)}
\def\kpa{k(\Phi,a)}
\def\LiE.{{\sf L\kern-.25em\raise0.59ex\hbox{\i}\kern-0.03em E}}

\def\calF{\mathcal{F}}
\def\calP{\mathcal{P}}
\def\calI{\mathcal{I}}
\def\cc{\mathfrak{c}}
\def\rd{R_\Delta}
\def\sd{S_\Delta}
\def\ais{{\alpha\in\sigma}}
\def\vol{\operatorname{vol}}

\def\al{\alpha}

\def\JKc{\operatorname{JK}_{\mathfrak{c}}}

\title{Vector partition function and representation theory}

\author[C. Cochet]{Charles Cochet}
\address{Institut de Math{\'e}matiques de Jussieu, 
         Universit{\'e} Paris 7 -- Denis Diderot, 
         Case 7012,
         75251 Paris cedex 05, France}
\email{cochet@math.jussieu.fr}
\date{January 1, 1994 and, in revised form, June 22, 1994.}


\begin{document}

\maketitle

\begin{abstract}
We apply some recent developments of 
Baldoni-Beck-Cochet-Vergne~\cite{BalBecCocVer05} on vector 
partition function, to Kostant's and Steinberg's formulae, 
for classical Lie algebras $A_r$, $B_r$, $C_r$, $D_r$.
We therefore get efficient {\tt Maple} programs that compute for 
these Lie algebras: 
the multiplicity of a weight in an irreducible finite-dimensional
representation; the decomposition coefficients of the tensor 
product of two irreducible finite-dimensional representations.
These programs can also calculate associated Ehrhart
quasipolynomials.
\medskip 

Nous appliquons des r{\'e}sultats r{\'e}cents de
Baldoni-Beck-Cochet-Vergne~\cite{BalBecCocVer05} sur la 
fonction de partition vectorielle, aux formules de 
Kostant et de Steinberg, dans le cas des alg{\`e}bres de Lie 
classiques $A_r$, $B_r$, $C_r$, $D_r$.
Ceci donne lieu {\`a} des programmes {\tt Maple} efficaces 
qui calculent pour ces alg{\`e}bres de Lie~:
la multiplicit{\'e} d'un poids dans une repr{\'e}sentation 
irr{\'e}ductible de dimension finie ; 
les coefficients de d{\'e}composition du produit tensoriel de 
deux repr{\'e}sentations irr{\'e}ductibles de dimension finie.
Ces programmes permettent {\'e}galement d'{\'e}valuer les
quasipolyn{\^o}mes d'Ehrhart associ{\'e}s.
\end{abstract}

%****************************************************
%****************************************************
\section{Introduction} 
\label{sect.intro}

In this note, we are interested in the two following computational
problems for classical Lie algebras $A_r$, $B_r$, $C_r$, $D_r$: 

\begin{itemize} 
\item The multiplicity $\clm$ of the weight $\mu$ in the 
  representation $V(\la)$ of highest weight $\la$. 
\item Littlewood-Richardson coefficients, that is the multiplicity
  $\clmn$ of the representation $V(\nu)$ in the tensor product of
  representations of highest weights $\la$ and $\mu$. 
\end{itemize} 

Softwares \LiE. (from van Leeuwen {\it et al.}~\cite{vanLee94})
and GAP~\cite{GAP}), and {\tt Maple} packages 
{\tt coxeter/weyl} (from Stembridge~\cite{Ste95}), use Freudenthal's
and Klimyk's formulae, and work for any semi-simple Lie algebra (not
only for classical Lie algebras).
Unfortunately, these formulae are really sensitive to the size of
coefficients of weights.
Moreover, they do not lead to the computation of associated
quasipolynomials $(\la,\mu)\mapsto\clm$ and
$(\la,\mu,\nu)\mapsto\clmn$.

Here the approach to these two problems is through 
\emph{vector partition function}, that is the function computing 
the number of ways one can decompose a vector as a linear 
combination with nonnegative integral coefficients of a fixed 
set of vectors.
For example the number $p(x)$ of ways of counting $x$ euros with 
coins, that is 
$$p(x)=\sharp\{n\in\N^8
   \,;\,
   x=200n_1+100n_2+50n_3+20n_4+10n_5+5n_6+2n_7+n_8\},$$

\noindent
can be seen as the partition of the $1$-dimensional vector $(x)$ 
with respects to the set 
$\{(200),(100),(50),(20),(10),(5),(2),(1)\}$ of 1-dimensional
vectors.
In the case of the decomposition with respects to the set of positive
roots of a simple Lie algebra, we speak of \emph{Kostant partition
  function}.

Recall that any $d$-dimensional rational convex polytope can be
written as the set $P(\Phi,a)$ of nonnegative solutions
$x=(x_i)\in\R^N$ of an equation $\sum_{i=1}^Nx_i\phi_i=a$, for a
matrix $\Phi$ with columns $\phi_i\in\Z^r$ and $a\in\Z^r$ ($d=N-r$).
It follows that evaluating the vector partition is equivalent to 
computing the number of integral points in a rational convex
polytope.

The vector partition function arises in many areas of mathematics: 
representation theory, flows in networks, magic squares, statistics, 
crystal bases of quantum groups.
Its complexity is polynomial in the size of input when the dimension 
of the polytope is fixed, and NP-hard if it can
vary~\cite{Bar94,Bar97,BarPom99}.

There are several approaches to the vector partition problem.
For example Barvinok's decomposition algorithm~\cite{Bar94}, 
recently implemented by the {\tt LattE}
team~\cite{DeLoeHemTauYos03,LattE}, works for general sets of
vectors.
Beck-Pixton~\cite{BecPix03} also created an algorithm dedicated 
to the vector set arising from the Birkhoff polytope, counting the
number of semi-magic squares.

In this note, we use recent results of 
Baldoni-Beck-Cochet-Vergne~\cite{BalBecCocVer05} to obtain 
a fast algorithm for Kostant partition function {\it via} 
inverse Laplace formula.
These results involve DeConcini-Procesi's
\emph{maximal nested sets} (or in short MNSs~\cite{DeConPro04})
and iterated residues of rational functions computed by formal
power series development.

We combine resulting procedures with Kostant's and Steinberg's
formulae giving $\clm$ and $\clmn$ in terms of vector partition
function.
We then obtain a {\tt Maple} program computing for classical Lie
algebras ($A_r$, $B_r$, $C_r$, $D_r$), the multiplicity of a weight
in an irreducible finite-dimensional representation,
as well as decomposition coefficients of the tensor product of
two irreducible finite-dimensional representations.
To the best of our knowledge, they are also the only ones able to
compute associated piecewise-defined quasipolynomials
$(\la,\mu)\mapsto\clm$ and $(\la,\mu,\nu)\mapsto\clmn$.

These programs (available at~\cite{home}) are specially designed
for large parameters of weights.
Indeed although only written in {\tt Maple} they can perform 
examples with weights with 5 digits coordinates,
far beyond classical softwares written in {\tt C++}.
We also stress that our programs are absolutely clear, easy to
use and require no installation of exotic package or program.
Retro-compatibility has been checked downto {\tt Maple Vr5}.
They are fully commented, so that a curious user can figure out
their internal mechanisms.

However, certain other softwares and packages are not limited by
the rank of the algebra like our programs.
For example computation of non-trivial examples in Lie algebras 
of rank $10$ is possible with the software \LiE., whereas 
our programs are efficient up to rank $5$--$7$.
These facts make our programs complementary to traditional 
softwares. 

Remark that Kostant's and Steinberg's formulae have already been 
implemented once in the case of $A_r$~\cite{Coc03a}.
This previous program relies on results of 
Baldoni-Vergne~\cite{BalVer01} implemented by 
Baldoni-DeLoera-Vergne~\cite{BalDeLoeVer03},
computing Kostant partition function only in the case of $A_r$.
Tools were \emph{special permutations} and again iterated residues
of rational fraction.

A new technique for Littlewood-Richardson coefficients has been
recently designed by DeLoera-McAllister~\cite{DeLoeMcAll05}.
For $A_r$, they wrote an algorithm using hive
polytopes~\cite{KnuTao99}.
For $B_r$, $C_r$, $D_r$, they implemented Berenstein-Zelevinsky
polytopes~\cite{BerZel01}.
They can also evaluate stretched Littlewood-Richardson
coefficients $\ctlmn$.
These two methods consist in computing a tensor product coefficient
as the number of lattice points in just one specific convex rational
polytope.
However our programs based on multidimensional residues are faster,
and can reach examples not available by their method.

This paper is organized as follows.
Section~\ref{sect.repr} recalls representation theory problems
we are interested in and links them with algebraic combinatorics.
Section~\ref{sect.poly} describes more precisely rational convex
polytopes and formulae counting their integral points.
Section~\ref{sect.mns} introduces maximal nested sets and
formulae that were used in our programs.
Finally in Section~\ref{sect.test} we perform tests of our programs.

%****************************************************
%****************************************************
\section{Representation theory and convex polytopes} 
\label{sect.repr}

Let us fix the notations once and for all.
Let $\g$ be a semi-simple Lie algebra of rank $r$.
Choose a Cartan subalgebra $\t$ of $\g$ and denote by $L\subset\t^*$
the weight lattice.

Let $\Delta^+$ be a positive roots system.
The root lattice is defined as $\Z[\Delta^+]$.
Let $C(\Delta^+)$ be the cone spanned by linear combinations with 
nonnegative coefficients of positive roots.
The Weyl group of $\g$ for $\t$ is denoted by $W$.

There exist only four simple Lie algebras $A_r$, $B_r$, $C_r$, $D_r$ of
rank $r$,
called \emph{classical Lie algebras} of rank $r$~\cite{Bou68},
and determined by their positive roots systems:

\begin{eqnarray*}
A_r:&&
  \Delta^+=\{e_i-e_j\,|\,1\leq i<j\leq r+1\}\subset\R^{r+1},\\
B_r:&&
  \Delta^+=\{e_i-e_j\,|\,1\leq i<j\leq r\}
             \cup\{e_i\,|\,1\leq i\leq r\}\subset\R^r,\\
C_r:&&
  \Delta^+=\{e_i-e_j\,|\,1\leq i<j\leq r\}
             \cup\{2e_i\,|\,1\leq i\leq r\}\subset\R^r,\\
D_r:&&
  \Delta^+=\{e_i-e_j\,|\,1\leq i<j\leq r\}
             \cup\{e_i+e_j\,|\,1\leq i<j\leq r\}\subset\R^r.
\end{eqnarray*}

The character of a representation $V$ of $\g$ is 
$\ch(V)=\sum_{\mu\in L}\dim(V_\mu)e^\mu$.
Recall that the irreducible finite-dimensional representation of $\g$ 
of highest weight $\la$ is denoted by $V(\la)$.
Hence the weight multiplicity $\clm$ is defined as $\dim(V(\la)_\mu)$
for any weight $\mu$ such that $\la-\mu$ is in the root lattice.
Multiplicities $\clm$ are called {\it Kostka numbers} when
$\g=A_r=\slr{r+1}$.

On the other hand, multiplicities of representations $V(\nu)$
in the tensor product $V(\la)\otimes V(\mu)$ are called
{\it Littlewood-Richardson coefficients}
(or {\it Clebsch-Gordan coefficients}).
Here $\nu$ is a dominant weight such that $\la+\mu-\nu$ is in the root
lattice.

Evaluating {\klr} is a difficult task.
For $A_1$, computing Kostka numbers is immediate and 
Clebsch-Gordan's formula gives Littlewood-Richardson coefficients.
For $A_2$, one can still compute some small examples. 
But for general $X_r$ ($r\geq 3$) or for weights which components 
are big (say, with two digits), direct computation is usually 
intractable.

There exist many formulae from representation theory for $\clm$ and 
$\clmn$. 
The first one, valid in any complex semi-simple Lie algebra $\g$, 
is Weyl's character formula

\begin{eqnarray*} 
\ch(V(\la)) & = & \frac{A_{\la+\rho}}{A_\rho}, 
\quad\mbox{ where } A_\mu=\sum_{w\in W}(-1)^{\eps(w)}e^{w(\mu)},
\end{eqnarray*} 

\noindent 
where $\rho$ is half the sum of positive roots for $\g$.
%Tensor product coefficients are obtained from this formula,
Littlewood-Richardson coefficients are obtained from this formula, 
since the character of $V(\la)\otimes V(\mu)$ is

$$\ch(V(\la)\otimes V(\mu)) \ = \ 
\ch(V(\la))\times\ch(V(\mu)) \ = \ 
  \sum_{\nu\in L\,;\,\la+\mu-\nu\in\Z[\Delta^+]}\clmn\ch(V(\nu)).$$

\noindent
But these two formulae do not lead to efficient computations when 
the rank of $\g$ or the size of coefficients of weights grow.
Moreover, computing the whole character is untractable:
for $\g=A_3=\slr 4$ and $\la=(2,1,0,-3)$, the character $\ch(V(\la))$
has 9 monomials but the character $\ch(V(10\la))$ has 2903
monomials.

Let us describe Kostant's and Steinberg's formulae in the case of 
any semi-simple Lie algebra $\g$.
Denote by $k_\g(a)$ the number of ways one can write a vector $a$ as
a nonnegative linear combination of positive roots.
Remark that $k_\g(a)=0$ unless $a$ is in the root lattice
$\Z[\Delta^+]$.
This number satisfies the equation 

\begin{eqnarray*}
\frac{1}{\prod_{\alpha\in\Delta^+}(1-e^{-\alpha})}
=\sum_{a\in\Z[\Delta^+]} k_\g(a)\ e^{-a}. 
\end{eqnarray*} 

Let $\la$ and $\mu$ be respectively a dominant weight and a weight
such that $\la-\mu\in\Z[\Delta^+]$.
A Weyl group element $w\in W$ is {\it valid} for $\la$ and $\mu$ if 
the root lattice element $w(\la+\rho)-(\mu+\rho)$ is in the cone
$C(\Delta^+)$. 
The set of such $w$'s is denoted by $\vlm$. 
Then Kostant's formula asserts that the weight multiplicity 
$\clm$ equals

\begin{equation} \label{equa.kost}
\clm=\sum_{w\in\vlm}(-1)^{\eps(w)}\,k_\g(w(\la+\rho)-(\mu+\rho)). 
\end{equation} 

\noindent
Similarly let $\la$, $\mu$, $\nu$, be three dominant weights such that
$\la+\mu-\nu\in\Z[\Delta^+]$.
The couple $(w,w')\in W\times W$ is {\it valid} for $\la$, $\mu$,
$\nu$, if the root lattice element
$w(\la+\rho)+w'(\mu+\rho)-(\nu+2\rho)$ is in $C(\Delta^+)$. 
The set of such couples is denoted by $\vlmn$.
Then Steinberg's formula asserts that the Littlewood-Richardson 
coefficient equals 

\begin{equation} \label{equa.stei}
  \clmn=\sum_{(w,w')\in\vlmn} (-1)^{\eps(w)+\eps(w')}
      k_{\g}(w(\la+\rho)+w'(\mu+\rho)-(\nu+2\rho)).
\end{equation} 

\noindent
Sets of valid Weyl group elements and valid couples of Weyl group 
elements turn out to be relatively small, when compared to $W$ and 
$W\times W$ (which size is exponential in the rank).
Remark that Kostant's (resp. Steinberg's) formula also work when
$\la-\mu$ (resp. $\la+\mu-\nu$) is not in the root lattice, since
Kostant partition function vanishes on vectors that are not in the
root lattice.

\medskip 

From now on, let $X_r$ be a classical Lie algebra of rank $r$.
Here $X$ stands for $A$, $B$, $C$, $D$.
Its positive roots system will be denoted by $X_r^+$.

\medskip 

Multiplicities $\clm$  and $\clmn$ behave nicely, in function of the
parameters.
More precisely, there exists a decomposition of the space
$\t^*\oplus\t^*\oplus\t^*$ in union of closed cones $C$, such that the
restriction of $\clmn$ to each cone $C$ is given by a
quasi-polynomial function.
This follows from theorems of Knutson-Tao~\cite{KnuTao99} (for $A_r$),
Berenstein-Zelevinsky~\cite{BerZel01} (for any semi-simple Lie
algebra) giving $\clmn$ as the number of points in a
rational convex polytope.
In the case of $A_r$, the fact that $\clmn$ is given on each cone $C$
by a polynomial function is proven in Rassart~\cite{Ras04}, and the
case of $A_3$ is treated as an illustration.
The description of the decomposition of $\t^*\oplus\t^*$ in cones
$C$, where the function $\clm$ is polynomial for $A_r$, was
given for low ranks by Billey-Guillemin-Rassart~\cite{BilGuiRas03}.
See also Rassart's website~\cite{Ras} for wonderful slides.

The common point to Kostant's and Steinberg's formulae is the function 
counting the number of decompositions of a root lattice element as a 
linear combination with nonnegative integral coefficients of
positive roots of the Lie algebra. 
The next section deals with an efficient method to compute it. 

%****************************************************
%****************************************************
\section{Counting integral points in rational convex polytopes} 
\label{sect.poly} 

\subsection{Vector partition function}

Let $E\simeq\R^r$ and $\Phi$ be an integral matrix with set of
columns $\Delta^+=\{\phi_1,\ldots,\phi_N\}\subset E^*$.
Choose $a\in\Z^r$.
The rational convex polyhedron associated to $\Phi$ and $a$ is

\begin{eqnarray*} 
\ppa & = & \left\{x\in\R^N\,;\,
  \sum_{i=1}^N x_i\phi_i=a,\,x_i\geq 0\right\}. 
\end{eqnarray*} 

\begin{remark} \label{rema.poly}
Every convex polyhedron can be realized under the form $\ppa$, 
that is as a set satisfying equality constraints on nonnegative 
variables.
Indeed any inequality can be replaced by an equality by adding
a new variable.
For example polytopes
$\{(x,y)\in\R^2\,;\,x\geq 0,y\geq 0,x+y\leq 1\}$ and
$\{(x,y,z)\in\R^3\,;\,x\geq 0,y\geq 0,z\geq 0,x+y+z=1\}$
are isomorphic and have the same number of integral points.
\end{remark}

We assume that $a$ is in the cone $C(\Phi)$ spanned by 
nonnegative linear combinations of the vectors $\phi_i$, so that 
$\ppa$ in non-empty. 
We also assume that the kernel of $\Phi$ intersects trivially
with the positive orthant $\R_+^N$, so that the cone 
$C(\Phi)$ is acute and $\ppa$ is a polytope ({\it i.e.} bounded). 
Finally, we assume that $\Phi$ has rank $r$.
The {\it vector partition function} is by definition

\begin{eqnarray*} 
k(\Phi,a) &  = & \left|P(\Phi,a)\cap\N^N\right|, 
\end{eqnarray*} 

\noindent 
that is the number of nonnegative integral solutions 
$(x_1,\ldots,x_N)$ of the equation $\sum_{i=1}^N x_i\phi_i=a$.
If $\Phi=\Phi(X_r)$ is the matrix which columns are positive roots
for a classical Lie algebra $X_r$, then 
$a\mapsto k(\Phi(X_r),a)$ is the {\it Kostant partition
function}.
For example 

$$\Phi(A_2)=
\left(\begin{array}{ccc}
 1 & 1 & 0\\
-1 & 0 & 1\\
 0 &-1 &-1
\end{array}\right)
\quad\mbox{and}\quad
\Phi(B_2)=
\left(\begin{array}{cccc}
 1 & 1 & 1 & 0\\
-1 & 1 & 0 & 1
\end{array}\right).$$

\noindent
Note that the matrix $\Phi(A_r)$ has rank $r$ (and not $r+1$), 
since sums on lines are zero.

\begin{figure}[ht]
\begin{center}
\includegraphics[height=3cm,width=4.5cm]{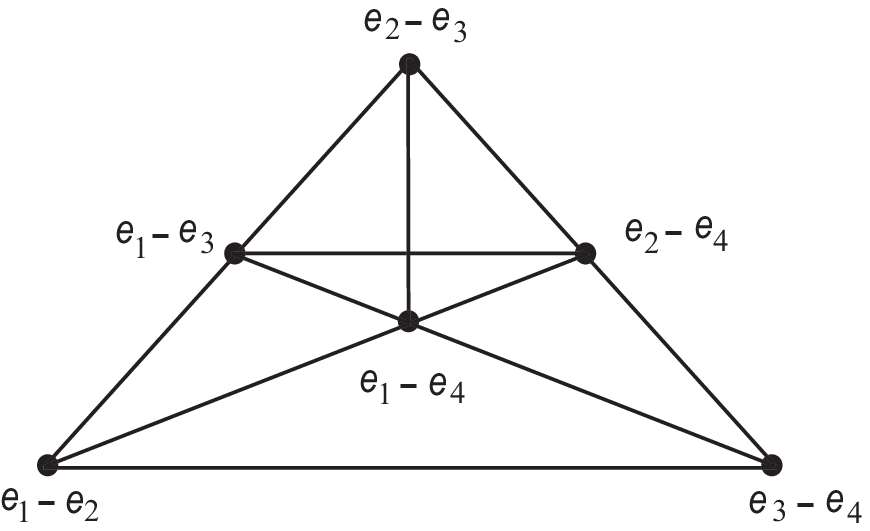}
\includegraphics[height=4cm,width=6cm]{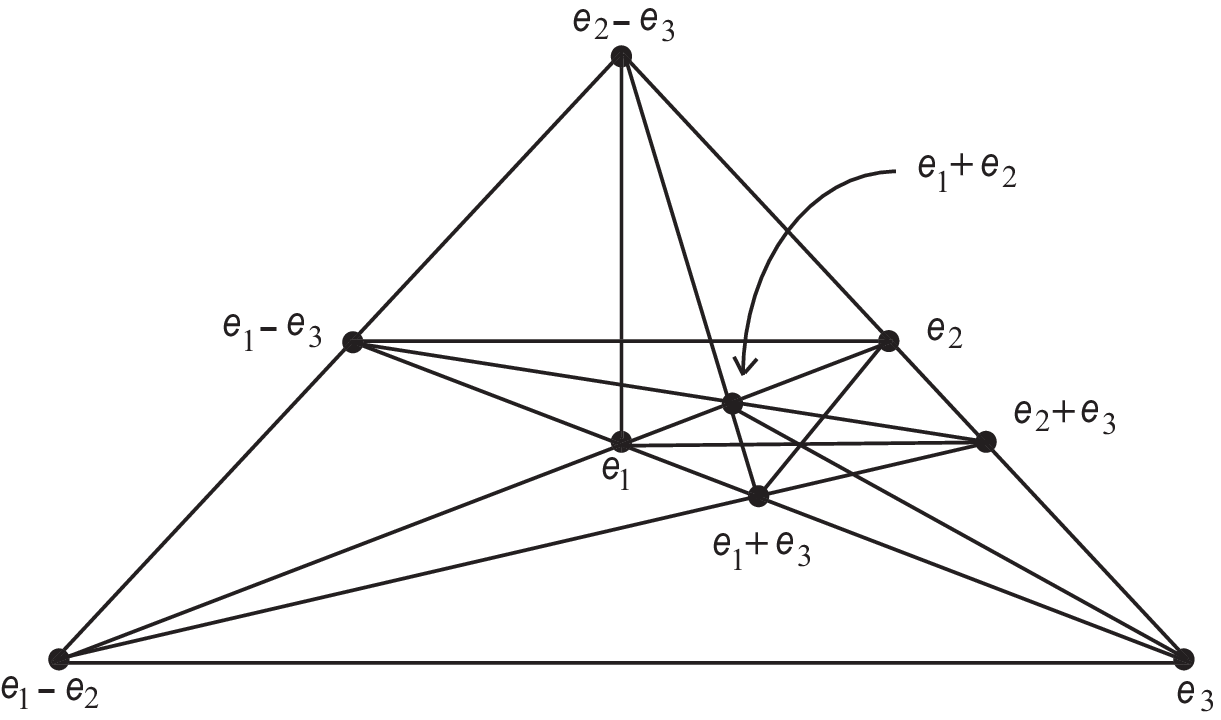}
\caption{The 7 chambers for $A_3$ and the 23 chambers for $B_3$}
\label{figu.A3B3}
\end{center}
\end{figure}

A \emph{basic subset} of $\Delta^+$ is a basis
$\sigma=\{\al_1,\ldots,\al_r\}$ of $E^*$ constituted with elements
of $\Delta^+$.
Let $B(\Delta^+)$ be the collection of all basic subsets of $\Delta^+$.
For such a $\sigma$, let $C(\sigma)$ be the cone of linear
combinations with nonnegative coefficients of $\al_i$'s.
Denote by $\operatorname{Sing}(\Delta^+)$ the reunion of the facets
of cones $C(\sigma)$, $\sigma\in B(\Delta^+)$; this is the set of
\emph{singular} vectors.
Let 
$C_{\operatorname{reg}}(\Delta^+)
  :=C(\Delta^+)\setminus\operatorname{Sing}(\Delta^+)$
be the set of \emph{regular} vectors.
A \emph{combinatorial chamber} $\cc$ is by definition a connected
component of $C_{\operatorname{reg}}(\Delta^+)$.
Combinatorial chambers are regions of quasi-polynomiality of
the vector partition function $a\mapsto\kpa$.
Figure~\ref{figu.A3B3} represents cones $C(A_3^+)$ and
$C(B_3^+)$, and their chamber decompositions.

%****************************************************
\subsection{Brion-Szenes-Vergne formula for classical Lie algebras} 

Let us describe the formula, computing the number of integral points
in rational convex polytopes $P(\Phi(X_r),a)$ associated to a
classical algebra $X_r$, that was implemented in our program.

Let $E=\t$ and consider the set $\Delta^+$ of positive roots
for $X_r$.
Denote by $\Delta$ the set $\Delta^+\cup(-\Delta^+)$ of all
roots.
Let $\rd$ be the vector space of fractions with poles on the
hyperplanes defined as kernels of forms $\alpha\in\Delta$.
Let $\sd$ be the vector space generated by fractions
$f_\sigma:=\frac{1}{\prod_{\ais}\alpha}$, $\sigma\in B(\Delta^+)$.
Brion-Vergne~\cite{BriVer97} proved that $\rd$ decomposes as the
direct sum $\sd\oplus\partial(\rd)$.
We define the \emph{Jeffrey-Kirwan residue} of the chamber $\cc$
as the linear form $\JKc$ on $\sd$:

$$\JKc(f_\sigma):=
\begin{cases}
  {\vol(\sigma)}^{-1}, &\mbox{if }\cc\subset C(\sigma),\\ 
  0,                   &\mbox{if }\cc\cap C(\sigma)=\emptyset,
\end{cases}
$$

\noindent
where $\vol(\sigma)$ is the volume of the parallelopiped
$\sum_{\alpha\in\sigma}[0,1]\alpha$.
We extend the JK residue to a linear form on $\rd$ by setting it
to $0$ on $\partial(\rd)$, and to a linear form on the space of
formal series $\widehat{\rd}$ by setting it to $0$ on
homogeneous elements of degree different from $-r$.
For example, for the system
$\Delta^+=\{e_1,e_2,e_1+e_2,e_1-e_2\}\subset\R^2$ of positive
roots for $B_2$ and the chamber $\cc=\N e_1\oplus\N(e_1+e_2)$
we have

$$
 \JKc\left(\frac{e^{x-y}}{xy^2}\right)
=\JKc\left(\frac{x-y}{xy^2}\right)
=\JKc\left(\frac{1}{y^2}-\frac{1}{xy}\right)
=-1,$$

\noindent
since $\cc\subset C(\{e_1,e_2\})$.

Let $T$ be the torus $E/E_\Z$, where $E_\Z\subset E$ is the dual of
the root lattice.
Given a basic subset $\sigma$, we define $T(\sigma)$ as the set of
elements $g\in T$ such that $e^{\langle \alpha,2i\pi G\rangle}=1$
for all $\alpha\in\sigma$;
here $G$ is a representative of $g\in E/E_\Z$.
Now let

$$\calF(g,a)(u):=
  \frac{e^{\langle a,2i\pi G+u\rangle}}
       {\prod_{\alpha\in\Delta}(1-e^{-\langle\alpha,2i\pi G+u\rangle})}.$$

\begin{theorem}[Brion-Szenes-Vergne~\cite{BriVer99,SzeVer04}]
  \label{theo.kpa}
Let $F\subset T$ be a finite set such that $T(\sigma)\subset F$ for all
$\sigma\in B(\Delta^+)$.
Fix a combinatorial chamber $\cc$.
Then for all $a\in\Z[\Delta^+]\cap\overline{\cc}$, we have:
$$\kpa=\sum_{g\in F}\JKc(\calF(g,a)).$$
\end{theorem}

Now that we linked vector partition function and Jeffrey-Kirwan
residue, we describe in Section~\ref{sect.mns} an efficient way to
compute the latter.

%****************************************************
%****************************************************
\section{DeConcini-Procesi's maximal nested sets (MNS)~\cite{DeConPro04}}
\label{sect.mns}

We keep the same notations as in Section~\ref{sect.poly}.
A subset $S\subset\Delta^+$ is \emph{complete} if
$S=\langle S\rangle\cap\Delta^+$.
A complete subset is \emph{reducible} if one can find a
decomposition $E=E_1\oplus E_2$ such that $S=S_1\cup S_2$ with
$S_1\subset E_1$ and $S_2\subset E_2$; else $S$ is said
\emph{irreducible}.
Let $\calI$ be the collection of irreducible subsets.

A collection $M=\{I_1,I_2,\ldots,I_s\}$ of irreducible subsets
$I_j$ of $\Delta^+$ is \emph{nested}, if:
for every subset $\{S_1,\ldots,S_m\}$ of $M$ such that there
exist no $i$, $j$ with $S_i\subset S_j$,
the union $S_1\cup S_2\cup\cdots\cup S_m$ is complete and the $S_i$'s
are its irreducible components.
Note that a maximal nested set (MNS in short) has exactly
$r$ elements.

Assume $\Delta^+$ irreductible and fix a total order on it.
For $M=\{I_1,\ldots,I_s\}$, $I_j\in\Delta^+$, take for every
$j$ the maximal element $\beta_j\in I_j$.
This defines an application
$\phi(M):=\{\beta_1,\ldots,\beta_s\}\subset\Delta^+$.
A maximal nested set $M$ is \emph{proper} if $\phi(M)$ is a basis
of $E^*$.
Denote by $\calP$ the collection of maximal proper nested sets
(MPNS in short).
We sort $\phi(M)$ and get an ordered list
$\theta(M)=[\alpha_1,\ldots,\alpha_r]$.
Thus $\theta$ is an application from the collection of MPNSs to
the collection of ordered basis of $E^*$.
For a given $M$, let then

\begin{eqnarray*}
   C(M)&:=&C(\alpha_1,\ldots,\alpha_r),\\
\vol(M)&:=&\vol\left(\oplus_{i=1}^r[0,1]\alpha_i\right),\\
\ires_M&:=&\res_{\alpha_r=0}\cdots\res_{\alpha_1=0}.
\end{eqnarray*}

\begin{example} \label{exam.Kn}
Let $e_i$ be the canonical basis of $\R^r$, with dual basis
$e^i$ ($i=1$, \ldots, $r$), and define $E$ as the subspace
of vectors which sum of coordinates vanish.
Consider the set $\Delta^+=\{e^i-e^j\,|\,1\leq i<j\leq r\}$ of
positive roots for $A_{r-1}$.
Irreducible subsets of $\Delta^+$ are indexed by subsets $S$ of
$\{1,2,\ldots,r\}$, the corresponding irreducible subset being
$\{e^i-e^j\,|\,i,j\in S, i<j\}$.
For instance $S=\{1,2,4\}$ parametrizes the set of roots given
by $\{e^1-e^2,e^2-e^4,e^1-e^4\}$.

A nested set is represented by a collection
$M=\{S_1,S_2,\ldots,S_k\}$ of subsets of $\{1,2,\ldots,r\}$ such
that if $S_i$, $S_j\in M$ then either $S_i\cap S_j$ is empty, or
one of them is contained in another.

For example one can easily compute that for the set of positive roots
for $A_3$ (see Figure~\ref{figu.A3B3}) there are only 7 MPNS, namely

$$\begin{array}{ll}
M_1=\{[1,2],[1,2,3],[1,2,3,4]\},&M_2=\{[2,3],[1,2,3],[1,2,3,4]\},\\
M_3=\{[2,3],[2,3,4],[1,2,3,4]\},&M_4=\{[3,4],[2,3,4],[1,2,3,4]\},\\
M_5=\{[1,3],[2,4],[1,2,3,4]\},  &M_6=\{[1,2],[3,4],[1,2,3,4]\}.
\end{array}$$
\end{example}

Now we can quote the Theorem for the Jeffrey-Kirwan residue computation:

\begin{theorem}[DeConcini-Procesi] \label{theo.jkc}
Let $\cc$ be a combinatorial chamber and fix $f\in\rd$.
Take any regular vector $v\in\cc$.
Then:
$$\JKc(f)
 =\sum_{M\in\calP\,:\,v\in C(M)}\frac{1}{\vol(M)}\ires_M(f).$$
\end{theorem}

See~\cite{BalBecCocVer05} for a detailed description of how
formulae from Theorems~\ref{theo.kpa} and~\ref{theo.jkc} were
implemented.

%****************************************************
%****************************************************
\section{Our programs}
\label{sect.test}

%****************************************************
\subsection{Description and implementation}

Initial data for weight multiplicity and
Littlewood-Richardson coefficients are only vectors
(respectively two and three).
Our programs work with weights represented in the canonical basis of
$E^*$, and not in the fundamental weights basis for $X_r$.
Translation between these two bases is performed {\it via}
straightforward procedures {\tt FromFundaToCanoX(r,v')} and
{\tt FromCanoToFundaX(r,v)} (where one replaces {\tt X} by {\tt A},
{\tt B}, {\tt C}, {\tt D}, according to the algebra).

Computation of the weight multiplicity $\clm$ and of the
Littlewood-Richardson coefficient $\clmn$ is done by typing in

\begin{eqnarray*} 
& & \mbox{{\tt MultiplicityX(lambda,mu)};}\\ 
& & \mbox{{\tt TensorProductX(lambda,mu,nu)};}\\ 
\end{eqnarray*} 

\noindent
where $\la$, $\mu$, $\nu$ are suitable weights.
The syntax for computing quasipolynomials is slightly
different.
Assume that we want to evaluate $(\la',\mu')\mapsto c_{\la'}^{\mu'}$
in a neighborhood of a couple $(\la,\mu)$,
and $(\la',\mu',\nu')\mapsto c_{\la'\,\,\,\,\mu'}^{\ \ \,\nu'}$ in a
neighborhood of a triple $(\la,\mu,\nu)$.
Let $\la_F=[x_1,\ldots,x_r]$, $\mu_F=[y_1,\ldots,y_r]$,
$\nu_F=[z_1,\ldots,z_r]$, be three formal vectors where $x_i$'s,
$y_i$'s and $z_i$'s are variables.
Then we use the command lines

\begin{eqnarray*}
& & {\tt PolynomialMultiplicityX(lambda,lambdaF,mu,muF);}\\
& & {\tt PolynomialTensorProductX(lambda,lambdaF,mu,muF,nu,nuF);}
\end{eqnarray*} 

\noindent
So for the polynomial $(\la',\mu')\mapsto c_{\la'}^{\mu'}$
with $\la=(3,2,1,-6)$ and $\mu=(2,2,-2,-2)$ for $A_3$ we enter

\begin{eqnarray*}
&&{\tt PolynomialMultiplicityA(}\\
&&\quad\quad{\tt [3,2,1,-6],[x[1],x[2],x[3],x[4]],[2,2,-2,-2],[y[1],y[2],y[3],y[4]]);}
\end{eqnarray*}

\noindent
and get instantly

$$\frac16(3x_1-2y_1+1)(3x_1-2y_1+2)(3x_1+6x_2-2y_1-6y_2+3).$$

\noindent
Remark that quasipolynomials $\ctlm$ and $\ctlmn$ are obtained
by setting $x_i=t\la_i$, $y_i=t\mu_i$, $z_i=t\nu_i$, so that

\begin{eqnarray*}
&&{\tt PolynomialMultiplicityA(}\\
&&\quad\quad{\tt [3,2,1,-6],[3t,2t,t,-6t],[2,2,-2,-2],[2t,2t,-2t,-2t]);}
\end{eqnarray*}

\noindent
returns $(t+1)(t+2)(t+3)/6$.

Now some words about implementation.
There are two main parts in our programs.
The first one is the implementation of Theorems~\ref{theo.kpa}
and~\ref{theo.jkc}; it is described in~\cite{BalBecCocVer05}.
The second one is the implementation of Kostant's~\eqref{equa.kost}
and Steinberg's~\eqref{equa.kost} formulae using valid Weyl group
elements and valid couples of Weyl group elements;
it is a generalization for classical Lie algebras of what has
been done for $A_r$ in~\cite{Coc03a}.

%****************************************************
\subsection{Comparative tests}

\begin{figure}[ht]
  \begin{center}
    \includegraphics[height=4cm,width=6cm]{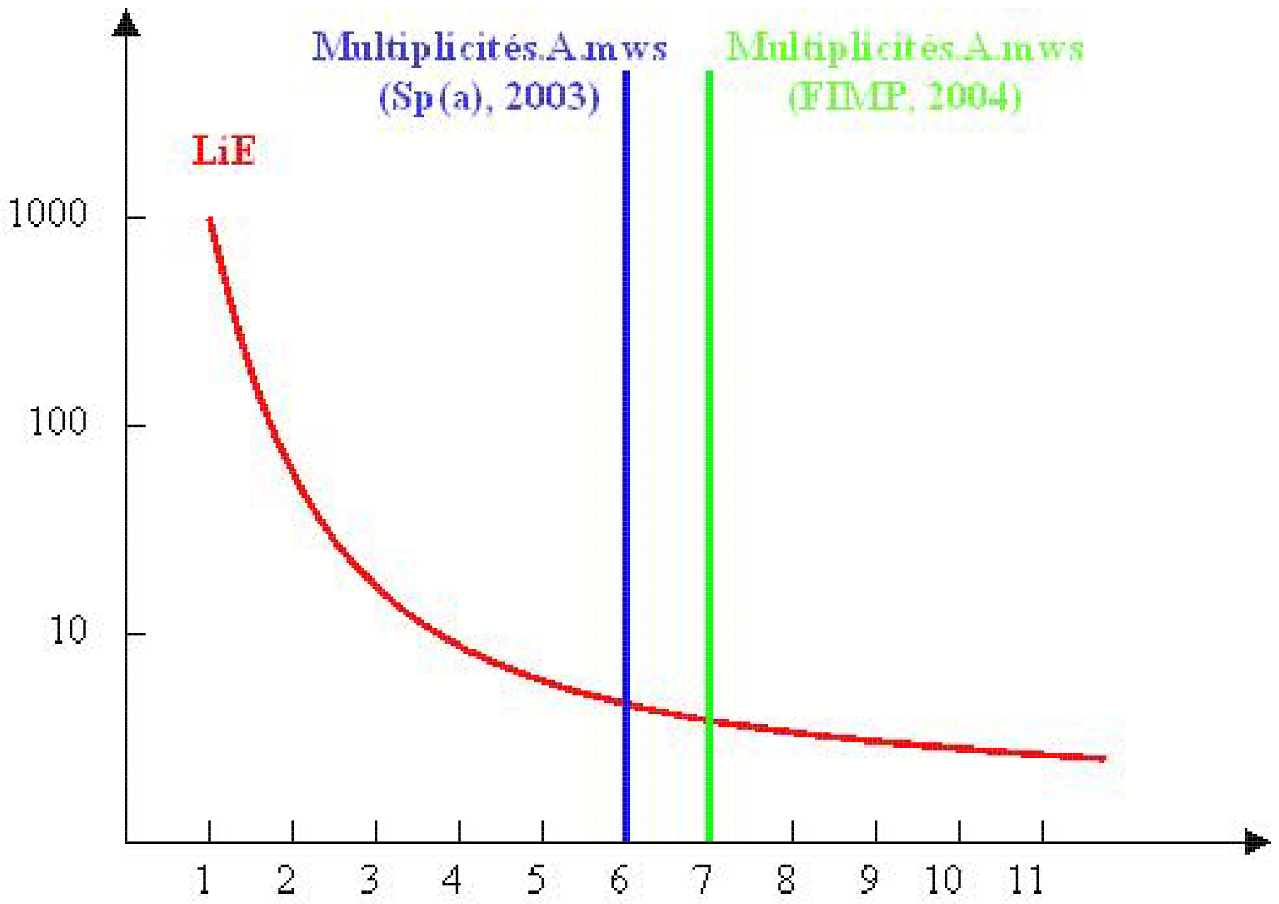}
    \includegraphics[height=4cm,width=6cm]{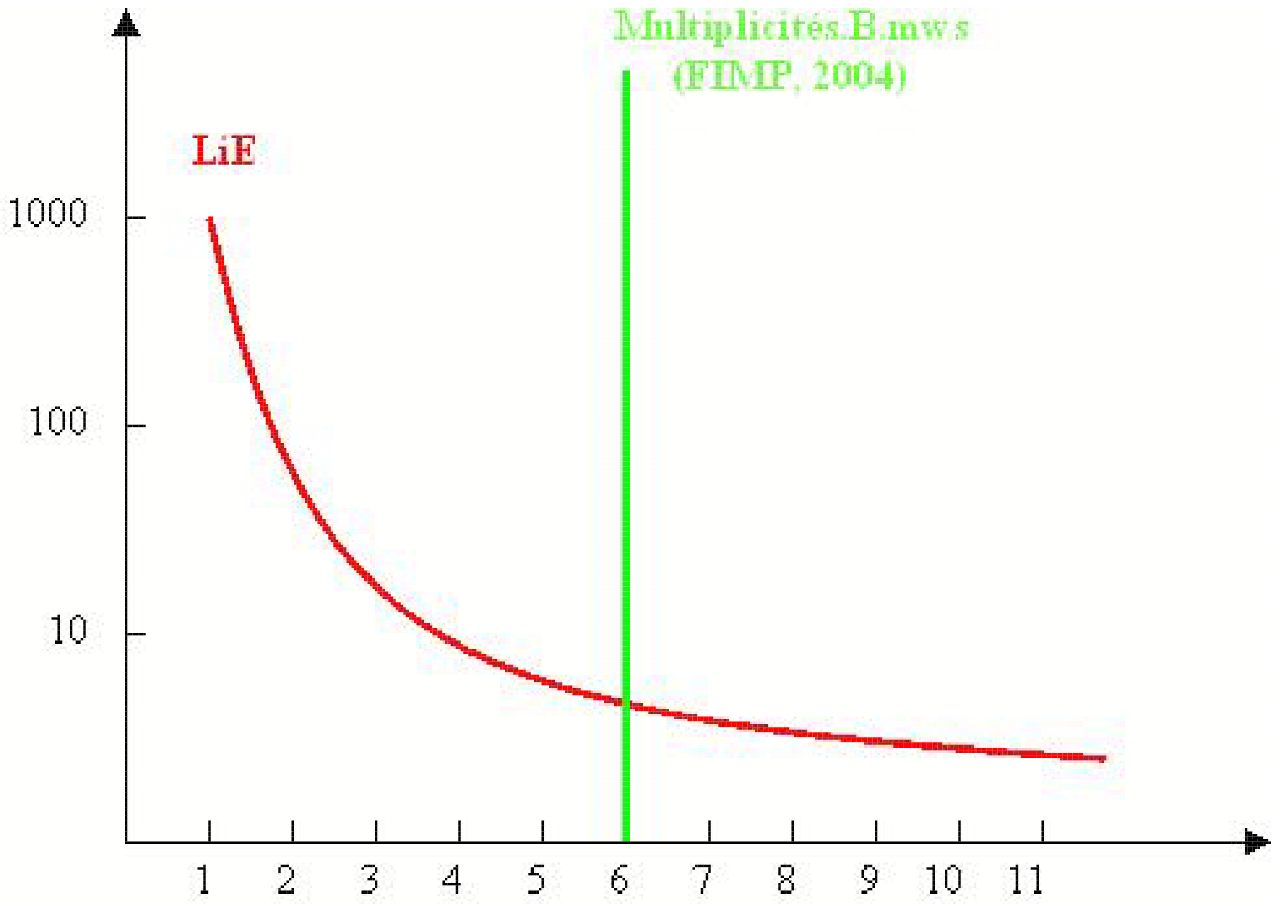}
    \caption{
      To the left, comparison for tensor product coefficients for $A_r$:
      with \LiE., with $\spa$ and with MNS.
      To the right, comparison for weight multiplicity of a weight for $B_r$:
      with \LiE. and with MNS.
      Similar Figures for $C_r$ and $D_r$.
    }
    \label{figu.comp}
  \end{center}
\end{figure}

Figure~\ref{figu.comp} describes efficiency area of the software
\LiE. and of our programs using MNS;
any area located to the left of a colored line represents the range 
where a program can compute examples in a reasonable time.
Figures~\ref{figu.A1}--\ref{figu.BCD} present precise comparative
tests of the software \LiE., of DeLoera-McAllister's
script~\cite{DeLoeMcAll05} using {\tt LattE}~\cite{LattE} and of
our programs using MNS.

All examples were runned on the same computer, a Pentium IV 1,13GHz
with 2Go of RAM memory.
Remark that computation times for {\tt LattE} and \LiE. are slower
than those shown in~\cite{DeLoeMcAll05}, due to different computers.
However, we performed exactly same examples for comparison purposes.

As in~\cite{DeLoeMcAll05}, in Tables~\ref{figu.A1}--\ref{figu.A2}
weights are for $\glr{r+1}$ and not $\slr{r+1}$ (coordinates
do not add to zero).
However the sum of coordinates of $\la+\mu-\nu$ vanish.

\begin{sidewaysfigure}
  \begin{tabular}{|l|r|r|r|r|}
    \hline
    $\la,\mu,\nu$ & $\clmn$ & {\tt MNS} & {\tt LattE} & {\tt LiE}\\
    \hline\hline
    (9,7,3,0,0), (9,9,3,2,0), (10,9,9,8,6)              &     2&8,0s& 3,0s&$<$0,1s\\
    (18,11,9,4,2), (20,17,9,4,0), (26,25,19,16,8)       &   453&2,8s& 8,8s&$<$0,1s\\
    (30,24,17,10,2), (27,23,13,8,2), (47,36,33,29,11)   &  5231&2,2s&11,4s&   0,5s\\
    (38,27,14,4,2), (35,26,16,11,2), (58,49,29,26,13)   & 16784&1,3s&12,8s&   1,5s\\
    (47,44,25,12,10), (40,34,25,15,8), (77,68,55,31,29) &  5449&1,3s& 8,8s&   1,4s\\
    (60,35,19,12,10), (60,54,27,25,3), (96,83,61,42,23) & 13637&1,0s& 8,4s&   9,1s\\
    (64,30,27,17,9), (55,48,32,12,4), (84,75,66,49,24)  & 49307&2,5s& 9,5s&  15,9s\\
    (73,58,41,21,4), (77,61,46,27,1), (124,117,71,52,45)&557744&2,1s&12,3s& 284,1s\\ \hline
  \end{tabular}
  \caption{For $A_r$, comparison of running times between the MNS
    algorithm, \texttt{LattE} and \LiE.}
  \label{figu.A1}
  
  \medskip
  \medskip

  \begin{tabular}{|l|r|r|r|}
    \hline
    $\la,\mu,\nu$ & $\clmn$ & {\tt MNS} & {\tt LattE}\\
    \hline
    \hline
    \begin{tabular}{l}
      (935,639,283,75,48)\\
      (921,683,386,136,21)\\
      (1529,1142,743,488,225)  % >25322,4s
    \end{tabular} 
    & 1303088213330 & 1,7s & 12,8s\\
    \hline
    \begin{tabular}{l}
      (6797,5843,4136,2770,707)\\
      (6071,5175,4035,1169,135)\\
      (10527,9398,8040,5803,3070)
    \end{tabular}
    & 459072901240524338 & 3,1s & 15,1s\\
    \hline
    \begin{tabular}{l}
      (859647,444276,283294,33686,24714)\\
      (482907,437967,280801,79229,26997)\\
      (1120207,699019,624861,351784,157647)
    \end{tabular}
    & 11711220003870071391294871475 & 2,0s & 11,9s\\
    \hline
  \end{tabular}
  \caption{For $A_r$, comparison of running times for large weights
    between the MNS algorithm and \texttt{LattE}}
  \label{figu.A2}
\end{sidewaysfigure}

\begin{sidewaysfigure}
  \begin{tabular}{|c|l|r|r|r|r|}
    \hline
         &$\la,\mu,\nu$&$\clmn$&{\tt MNS}&{\tt LattE} &{\tt LiE}\\ \hline \hline
    $B_3$&(46,42,38), (38,36,42), (41,36,44)        &   354440672&  6,4s& 22,5s&  229,0s\\
         &(46,42,41), (14,58,17), (50,54,38)        &    88429965&  2,7s& 15,2s&  102,6s\\
         &(15,60,67), (58,70,52), (57,38,63)        &   626863031&  7,8s& 17,0s&  713,5s\\
         &(5567,2146,6241), (6932,1819,8227), (3538,4733,3648)&215676881876569849679&7,0s&16,3s&--\\ \hline
    $C_3$&(25,42,22), (36,38,50), (31,33,48)        &    87348857&  5,6s& 18,1s&   52,9s\\
         &(34,56,36), (44,51,49), (37,51,54)        &   606746767&  5,1s& 20,4s&  516,0s\\
         &(39,64,58), (65,15,72), (70,41,44)        &   519379044&  8,7s& 18,3s& 1096,9s\\
         &(5046,5267,7266), (7091,3228,9528), (9655,7698,2728)&1578943284716032240384&8,2s&18,3s&--\\ \hline
    $D_4$&(13,20,10,14), (10,20,13,20), (5,11,15,18)&    41336415&131,0s&185,8s&  224,7s\\
         &(12,22,9,30), (28,14,15,26), (10,24,10,26)&   322610723& 78,6s&192,7s& 1184,8s\\
         &(37,16,31,29), (40,18,35,41), (36,27,19,37)&18538329184& 64,3s&258,7s&21978,4s\\
         &(2883,8198,3874,5423), (1901,9609,889,4288), (5284,9031,2959,5527)&1891293256704574356565149344&27,7s&165,2s&--\\
         \hline
  \end{tabular}
  \caption{For $B_r$, $C_r$, $D_r$, comparison of running times
    between {\tt LattE}, the MNS algorithm and \LiE.}
  \label{figu.BCD}
\end{sidewaysfigure}

Now some words about quasipolynomials computation.
Let us examine the first example for $B_3$ in~\cite{DeLoeMcAll05},
that is the evaluation of the quasipolynomial $\ctlmn$ for weights
$\la=[0,15,5]$, $\mu=[12,15,3]$ and $\nu=[6,15,6]$ expressed in the
basis of fundamental weights.
In canonical basis, these data become
$\la=(35/2,35/2,5/2)$, $\mu=(57/2,33/2,3/2)$, $\nu=(24,18,3)$.
The program using the MNS algorithm returns the quasipolynomial

\begin{eqnarray*}
\ctlmn
&=& \left(\frac{203}{256}+\frac{53}{256}(-1)^t\right)
   +\left(\frac{1515}{128}+\frac{197}{128}(-1)^t\right)t\\
& &+\left(\frac{35353}{384}+\frac{881}{128}(-1)^t\right)t^2
   +\left(\frac{13405}{32}\right)t^3\\
& &+\left(\frac{407513}{384}\right)t^4
   +\left(\frac{68339}{64}\right)t^5
\end{eqnarray*}

\noindent
in 1099,4s.
On the other hand, the computation of the full quasipolynomial $\clmn$
with formal vectors $[x_1,x_2,x_3]$, $[y_1,y_2,y_3]$, $[z_1,z_2,z_3]$
leads to a 87 pages result, obtained in only 1158,6s.
With {\tt LattE}, on our computer, one obtains the
quasipolynomial $\ctlmn$ in only 825,8s.

As announced in the introduction, our program is really efficient for
weights with huge coefficients.
Note that in the particular case of $A_r$ the MNS algorithm
allows us to compute examples one rank further than the $\spa$
algorithm.

The translation of the program using MNS in the language of the
symbolic calculation software {\tt MuPAD} is in progress.
A version using distributed calculation on a grid of computers is in
the air; it will considerably increase the speed of computations.

%****************************************************
%****************************************************

\bibliographystyle{amsalpha}

\end{document}